\def\fnote#1{\footnote}
\documentstyle[12pt,leqno]{article}
\newtheorem{th}{Theorem}[section]
\newtheorem{lem}[th]{Lemma}
\newtheorem{pro}[th]{Proposition}
\newtheorem{co}[th]{Corollary}
\newtheorem{rem}[th]{Remark}

\title{\Large \bf {Analytic fields on compact balanced\\
Hermitian manifolds}}
\author{{\sc G. Ganchev}  \thanks{The author
supported by Contract
MM 413/1994 with the Ministry of Science and Education of
Bulgaria.}
\hspace{3mm} {\sc S.Ivanov} \thanks{The author supported
by Contract
MM 413/1994 with the Ministry of Science and Education of
Bulgaria and by
Contract 219/1994 with the University of Sofia "St. Kl.
Ochridski".}}
\date{}
\begin{document}
\maketitle
\thispagestyle{empty}
\vspace{25mm}
\begin{abstract}
On a Hermitian manifold we construct a symmetric $(1,1)$-
tensor $H$ using the torsion and the curvature of the Chern
connection. On a compact balanced Hermitian
manifold we find necessary and
sufficient conditions in terms of the tensor $H$ for a
harmonic $1$-form to be analytic and for an analytic $1$-
form to be
harmonic. We prove that if $H$ is positive
definite then the first Betti number $b_1 = 0$
and the Hodge number $h^{1,0} = 0$. We obtain an obstruction 
to the
existence of Killing vector
fields in terms of the Ricci tensor  of the
Chern connection.
We prove that if the Chern form of the
Chern connection on
a compact balanced Hermitian manifold is non-positive 
definite then
every Killing
vector field is analytic; if moreover the Chern form is
negative definite then there
are no Killing vector fields.
\\[15mm]
{\bf Running title:} Analytic fields on balanced
manifolds\\[5mm]
{\bf Keywords.} Compact balanced manifolds,
Killing, affine, holomorphic vector
fields, harmonic, holomorphic $1$-forms, vanishing theorem
of Bochner type.
\\[5mm]
${\bf MS}$ {\bf classification: } 53C15; 53C55; 53B35
\end{abstract}
\newpage
\section{Introduction}
On a compact Riemannian manifold $({\bf M},g)$ certain
objects of geometric
interest, such as Killing vector fields and harmonic $1$-forms
must satisfy
additional differential equations when appropriate Levi-
Civita curvature
conditions are imposed. This leads to obstructions to
the existence of
these objects, known as  Bochner-type vanishing theorems.
For example, the
well known theorems of Bochner state that if the Ricci
tensor
is non-negative
(resp. non-positive) then every harmonic 1-form (resp. every
Killing vector
field) must be parallel and if the Ricci tensor is negative
(resp. positive) at
one point then there are no harmonic $1$-forms (resp. no
Killing vector
fields). By the Hodge theory the nonexistence of harmonic
$1$-forms leads to the vanishing of
the first de Rham
cohomology group $H^1({\bf M},{\bf R})$.
\par
On a Hermitian manifold $({\bf M},g,J)$ another objects of
geometric interest are holomorphic vector fields and
holomorphic forms. Closely
related to the hermitian structure is the Chern connection.
On a compact
Hermitian manifold holomorphic vector fields and holomorphic
$(p,0)$-forms
have
to satisfy additional differential equations when
appropriate curvature
conditions on the Chern connection are imposed. This also
leads to some
obstructions to the existence of these objects,  so
called
vanishing theorems
for holomorphic sections. For example, if the mean curvature
of the Chern
connection (in the sense of \cite{Ko1}) on a compact
Hermitian manifold is
non-negative (resp. non-positive) then every holomorphic
$(1,0)$-form (resp.
every holomorphic vector field) is parallel with respect to
the Chern
connection and if the mean curvature is positive (resp.
negative) at one point
then there are no holomorphic $(1,0)$-forms (resp. no
holomorphic vector
fields) (for more general formulations see \cite{KW,Ga}). By 
the
Hodge theory (see
\cite{H}) this leads
to the vanishing of the Dolbeault cohomology group
$H^{1,0}({\bf M},{\bf C})$.
\par
In the Kaehler geometry there is  a close relation  between
the metric objects
(Killing vector fields, harmonic $1$-forms) on one hand and
the holomorphic objects
(holomorphic vector fields, holomorphic $(1,0)$-forms) on
the other hand.
For a Kaehler manifold the Levi-Civita connection coincides
with the Chern
connection. Then the mean curvature of the Chern connection
is exactly the
Ricci tensor. Thus, on a compact Kaehler manifold the
positive (resp. negative)
definitness of the Ricci tensor is an obstruction to the
existence not
only of harmonic $1$-forms (resp. Killing vector fields)
but also to the
existence of holomorphic $(1,0)$-forms (resp. holomorphic
vector fields).
Moreover, on a compact Kaehler manifold an $1$-form is
harmonic iff it is
analytic (i.e its $(1,0)$-part is holomorphic)  and every
Killing vector field is analytic (see \cite{Ya2}). So the
first Betti mumber
$b_1 = dim H^1({\bf M},{\bf R})$ vanishes iff the Hodge
number $h^{1,0} = dim
H^{1,0}({\bf M},{\bf C})$ vanishes.
\par
In general, on a compact Hermitian manifold there is no
remarkable relation
between Killing and holomorphic vector fields and between
harmonic $1$-forms
and holomorphic $(1,0)$-forms. \par
In this paper we consider compact balanced Hermitian
manifolds and try to find
a connection between the metric objects and the holomorphic
objects mentioned above.
Balanced manifolds have been studied intensively in
\cite{Mi,A_B1,A_B2,A_B3};
in \cite{Ga} they are called semi-Kaehler of special type.
This class of manifolds
includes the class of Kaehler manifolds but also many
important classes of
non-Kaehler manifolds, such as: complex solvmanifolds,
twistor spaces of oriented riemannian 4-manifolds, 1-
dimensional families of Kaehler manifolds
(see \cite{Mi}), hermitian compact manifolds with flat Chern
connection (see
\cite{Ga}), twistor spaces of oriented distinguished Weyl
structure on compact self-dual 4-manifolds \cite{Ga1},
twistor spaces of quaternionic Kaehler manifolds
\cite{P,AGI}, manifolds obtained as modification of compact
Kaehler manifolds
\cite{A_B1} and of compact balanced manifolds \cite{A_B2},
see also
\cite{A_B3}.
\par
We construct on a Hermitian manifold  a symmetric $(1,1)$-
tensor $H$ using the torsion and the curvature of the Chern
connection. On a compact balanced Hermitian
manifold we give in Theorem \ref{th42} necessary and
sufficient conditions in terms of the tensor $H$ for a
harmonic $1$-form to be analytic and for an analytic $1$-
form
to be harmonic.
This alows us to obtain a vanishing theorem of Bochner type
on compact balanced Hermitian manifolds
(Theorem \ref{th46}). We prove that if $H$ is positive
definite then $b_1 = 0$
and $h^{1,0} = 0$. We obtain an obstruction to the
existence of Killing vector
fields in terms of the Ricci tensor (or the Chern form) of the
Chern connection. In
Theorem \ref{th39} we prove that if the Chern form of the
Chern connection on
a compact balanced Hermitian manifold is non-positive then
every Killing
vector field is analytic; if moreover the Chern form is
negative then there
are no Killing vector fields.\par
It is well known that on a compact Riemannian manifold a
smooth vector field is
Killing iff it is affine with respect to the Levi-Civita
connection {\cite{Ya2}. Thus, on a compact Kaehler manifold
every affine vector field is analytic.  On a compact
balanced Hermitian manifold we find
necessary and sufficient conditions in terms of the Lie
derivative of the Chern
connection for a smooth vector field to be analytic. In
particular we prove that every affine vector field with
respect to the Chern
connection is analytic
on a compact balanced Hermitian manifold.
\section{Preliminaries}
Let $({\bf M},J,g)$ be a $2n$-dimensional Hermitian manifold
with complex
structure $J$ and
Riemannian metric $g$.The algebra of all $C^{\infty }$
vector fields on {\bf
M}  will be denoted by  {\bf XM}. The complex structure $J$
on the tangent
bundle ${\bf TM}$ of ${\bf M}$ induces a splitting of the
complexified tangent
bundle ${\bf T}_c{\bf M}$ into two complementary subbundles,
conjugate to each
other:${\bf T}_c{\bf M} = {\bf T}^{(1,0)}{\bf M} + {\bf
T}^{(0,1)}{\bf M}$. The elements of ${\bf T}^{(1,0)}{\bf M}$
(resp. ${\bf T}^{(0,1)}{\bf M}$) are the (complex)
tangent vectors of type $
(1,0)$ (resp. of type $(0,1)$).Each real tangent vector
field $X$ can be
expressed  in a unique way  as a sum: $X = U + \bar{U}$,
where $U =
\frac{1}{2} (X - \sqrt {-1}JX) \in {\bf T}^{(1,0)}{\bf M}$
and
$\bar{U} = \frac{1}{2} (X +
\sqrt {-1}JX) \in {\bf T}^{(0,1)}{\bf M}$. With respect to
local
holomorphic coordinates
$\{z^{\alpha }\}$, ($\alpha = 1,...,n$) we have $U =
X^{\alpha } \frac{\partial
}{\partial z^{\alpha }}, \bar{U} = X^{\bar{\alpha }}
\frac{\partial }{\partial
Z^{\bar{\alpha }}}$ (summation
convention is assumed further in the paper). The induced
complex structure  on
the cotangent bundle ${\bf T}^{*}{\bf M}$ (also denoted by
$J$) is defined by:
$(J\omega )(X) = - \omega (JX)$, where $\omega $ is a real
1-form and $X$ is a
real vector field on ${\bf M}$. For the complexified
cotangent bundle ${\bf
T}^{*}_c{\bf M}$ we have the splitting: ${\bf T}^{*}_c{\bf
M} = \Lambda ^{1,0}({\bf
M}) + \Lambda ^{0,1}({\bf M})$. The elements of $\Lambda
^{1,0}({\bf M})$ (resp.
$\Lambda ^{0,1}({\bf M})$) are the (complex) 1-forms of type
$(1,0)$ (resp. of
type $(0,1)$. Each real 1-form $\omega $ can be expressed
in a unique way
as a sum $\omega = \beta + \bar{\beta }$, where $\beta =
\frac{1}{2}(\omega -
\sqrt {-1}J\omega ) \in \Lambda ^{1,0}({\bf M}); \bar{\beta
} = \frac{1}{2}(\omega +
\sqrt {-1}J\omega ) \in \Lambda ^{0,1}({\bf M})$. With
respect to  local
holomorphic coordinates we have $\beta = \omega _{\alpha
}dz^{\alpha },
\bar{\beta } = \omega _{\bar{\alpha }}dz^{\bar{\alpha }}$.
In the whole paper
all  tensors and connections will be extended complex
multilinearly to the
complexification  ${\bf T}_c{\bf M}$ of ${\bf T}{\bf M}$.
\par
The Kaehler form $\Omega $ on {\bf M} is
defined by $\Omega (X,Y) = g(JX,Y); X,Y \in {\bf XM}$. The
Lee form $\theta $ of
the Hermitian structure is defined by $\theta = - \delta
\Omega \circ J$.
\par
The Levi-Civita connection and the canonical Chern
connection (the Hermitian
connection) will be denoted by $\nabla $ and $D$,
respectively. We recall that
the Chern connection $D$ is the unique linear connection
preserving the metric
and the complex structure with torsion tensor $T$, having
the following
property: $T(JX,Y) = T(X,JY)$. This  implies (e.g.
\cite{A-Z}):
\begin{equation}\label{1}
T(JX,Y) = JT(X,Y), \qquad X,Y \in {\bf XM}.
\end{equation}
The two connections are related by the following identity
\begin{equation}\label{2}
g(\nabla _XY,Z) = g(D_XY,Z) + \frac{1}{2} d\Omega (JX,Y,Z),
\qquad X,Y,Z \in
{\bf XM}.
\end{equation}
Let $e_{1},...,e_{2n}$ be an orthonormal local basis on {\bf
M}. We consider the
following Ricci-type tensors associated with the curvature
tensor $K$ of the
Chern connection:
$$k(X,Y) = - \frac{1}{2}\sum
_{j=1}^{2n}g(K(X,JY)e_{j},Je_{j});
\quad k^{*}(X,Y) = -
\frac{1}{2}\sum _{j=1}^{2n}g(K(e_{j},Je_{j})X,JY);$$
$$ s(X,Y) = \sum_{j=1}^{2n}g(K(e_{j},X)Y,e_{j}).$$
The $(1,1)$-form $\kappa $ corresponding to the tensor $k$
represents the first Chern
class of ${\bf M}$ (further we will call it the Chern form)
and the (1,1)-form $\kappa ^*$
corresponding to the tensor $k^*$ is the mean curvature of
the holomorphic
tangent bundle ${\bf T}^{(1,0)}{\bf M}$ with the  hermitian
metric
$g$.
\par
Using the torsion tensor $T$ of the Chern connection we
construct another
remarkable symmetric (1,1) tensor as follows
$$
t(X,Y) = \sum_{\alpha ,\beta =1}^{n}g(T(E_{\alpha },E_{\beta
}),X)g(T(E_{\bar{\alpha }},E_{\bar{\beta}}),Y),
$$
where $E_1,...,E_n,E_{\bar{1}},...,E_{\bar{n}}$ is a
hermitian basis on ${\bf
T}_c{\bf M}$. From the definition and (\ref{1}) it follows
that $t$ is
symmetric, $J$-
invariant positive semi-definite tensor i.e. $t(X,Y) =
t(Y,X) = t(JX,JY);
\quad t(X,X) \ge 0, X,Y \in {\bf XM}$
\par
As we shall see below, the tensor $t$ plays an important
role for the
coinsidence of harmonic 1-forms with analytic 1-forms. \par
For a real 1-form $\omega $ we denote by $\omega ^{\#}$ 
the
corresponding
vector fiel defined by:$\omega (Y) := g(\omega ^{\#},Y), Y
\in {\bf XM}$ and
for
a real vector field $X$ we denote by $\omega _X$ the
corresponding real 1-form
defined by: $\omega _X(Y) := g(X,Y)$. If $\omega = \omega
_{\alpha
}dz^{\alpha }$ is an $(1,0)$-form (resp. $\omega
_{\bar{\alpha
}}dz^{\bar{\alpha
}}$ is a $(0,1)$-form) then $\omega ^{\#} = g^{\alpha
\bar{\beta }}\omega
_{\alpha
}\frac{\partial }{\partial z^{\bar{\beta }}}$ is the
corresponding
$(0,1)$-vector field
(resp. $\omega ^{\#} = g^{\alpha \bar{\beta }}\omega
_{\bar{\beta
}}\frac{\partial
}{\partial z^{\alpha }}$ is the $(1,0)$-vector field) and if
$X = X^{\alpha
}\frac{\partial }{\partial z^{\alpha }}$ is an $(1,0)$-
vector field (resp. $X
=
X^{\bar{\alpha }}\frac{\partial }{\partial z^{\bar{\alpha
}}}$ is a
$(0,1)$-vector
field) then $\omega _X = g_{\alpha \bar{\beta }}X^{\alpha
}dz^{\bar{\beta }}$
is a $(0,1)$-form (resp. $\omega _X = g_{\bar{\alpha }\beta
}X^{\bar{\alpha
}}dz^{\beta }$ is an $(1,0)$-form).
\par
For a real 1-form $\omega $ using (\ref{2}) we calculate
\begin{equation}\label{3}
(\nabla _X\omega )Y = (D_X\omega)Y - \frac{1}{2}d\Omega
(JX,Y,\omega ^{\#}),
\qquad X,Y \in {\bf XM}.
\end{equation}
>From (\ref{3}) it follows that
\begin{equation}\label{4}
\delta \omega = - \sum_{i=1}^{2n}(D_{e_i}\omega )e_i -
\theta (\omega
^{\#}). \end{equation}
\indent
We recall here the definition of a
balanced manifold and some characterizations given in
\cite{Mi} and \cite{Ga}, for
completeness: \par
DEFINITION: {\it A balanced} manifold {\bf M} is a compact
complex $n$-
manifold which satisfies one of the following equivalent
conditions:\par
i) {\bf M} admits a hermitian structure $(g,J)$ such that
$d\Omega ^{n-1} = 0$;
\par
ii) {\bf M} admits a hermitian structure $(g,J)$ such that
$\delta \Omega = 0$;
\par
iii) {\bf M} admits a hermitian structure $(g,J)$ such that
$\theta = 0$; \par
iv) {\bf M} admits a hermitian structure $(g,J)$ such that
$\Delta_{\partial}f =
\Delta_{\bar{\partial}}f = \frac{1}{2}\Delta_df$, for any
smooth function $f$
on
${\bf M}$, where $\Delta_{\partial},
\Delta_{\bar{\partial}}, \Delta_df$
denote
the Laplacians with respect to the operators $\partial,
\bar{\partial},d$,
respectively.
\par
v) there are no non-zero positive $(n-1,n-1)$-currents on
${\bf M}$ wich are $(n-
1,n-1)$-components of boundaries. \par
In this paper we shall use essentially iii) and iv).\par
On balanced manifolds the first Ricci tensor $k$ coincides
with the third Ricci
tensor $s$ (see e.g. \cite{Ba}).\par
If $({\bf M},J,g)$ is a balanced manifold then the equality
(\ref{4}) takes
the form \begin{equation}\label{5}
\delta \omega = - \sum_{i=1}^{2n}(D_{e_i}\omega )e_i .
\end{equation}
\section{Analytic and Killing vector fields}

A real vector field $\xi $ is said to be {\it analytic} if
$L_{\xi}J = 0$,where
$L_{\xi }$ denotes the Lie derivative with respect to $\xi
$. The vector field
$\xi $ is analytic iff the $(2,0)$ part of $D\omega _{\xi }$
vanishies. In local
holomorphic coordinates this condition can be written as
follows
\begin{equation}\label{6}
D_{\alpha }\xi _{\beta } = 0,
\end{equation}
i.e. the $(1,0)$-part of $\xi $ is a holomorphic vector field.
We consider the following real 1-form $\omega $ defined by
$$
\omega = \xi ^{\beta}D_{\alpha }\xi _{\beta }dz^{\alpha } +
\xi ^{\bar{\beta
}}D_{\bar{\alpha }}\xi _{\bar{\beta }}dz^{\bar{\alpha }}.
$$
Using essentially that $DJ = 0$ and (\ref{3}), we find
$$
\delta \omega = - \Vert D_{\alpha }\xi _{\beta } \Vert^2 -
2Re\left[ \xi ^{\beta
}D^{\alpha }D_{\alpha }\xi _{\beta } + \xi ^{\beta }\theta
^{\alpha }D_{\alpha
}\xi _{\beta } \right],
$$ where $2Re(f)$ denotes the real part of a complex valued
function $f$ and $\Vert D_{\alpha }\xi _{\beta } \Vert^2$ is
the norm of the $(2,0)+(0,2)$-part of $D\omega _{\xi }$. The
norm of the $(1,1)$-part of $D\omega _{\xi }$ will be
denoted by $\Vert D_{\alpha }\xi _{\bar{\beta }} \Vert^2$.\\
Thus, on a compact Hermitian manifold we have the formula
\begin{equation}\label{7}
\int_{{\bf M}} \left\{ \Vert D_{\alpha }\xi _{\beta }
\Vert^2 + 2Re\left[ \xi
^{\beta
}D^{\alpha }D_{\alpha }\xi _{\beta } + \xi ^{\beta }\theta
^{\alpha }D_{\alpha
}\xi _{\beta } \right]\right\}  \,dV = 0.
\end{equation}
Taking into account the Ricci formula
\begin{equation}\label{7*}
D^{\alpha }D_{\alpha }\xi _{\beta } = D_{\alpha }D^{\alpha
}\xi _{\beta } +
k^*_{\beta \bar{\sigma }}\xi ^{\bar{\sigma }}
\end{equation}
from (\ref{7}) we obtain
\begin{pro}\label{pro31}
Let $\xi $ be a real vector field on a compact Hermitian
manifold. The following
conditions are equivalent:
\par i) $\xi $ is analytic;
\par ii) $D^{\alpha }D_{\alpha }\xi _{\beta } + \theta
^{\alpha }D_{\alpha
}\xi _{\beta } = 0;$
\par iii) $D_{\alpha }D^{\alpha }\xi _{\beta } + k^*_{\beta
\bar{\sigma }}\xi
^{\bar{\sigma }} + \theta ^{\alpha }D_{\alpha }\xi _{\beta }
= 0.$
\end{pro}
Hence, on a compact balanced Hermitian manifold we have
\cite{A-Z}
\begin{pro}\label{pro32}
Let $\xi $ be a real vector field on a compact balanced 
Hermitian
manifold. The following
conditions are equivalent:
\par i) $\xi $ is analytic;
\par ii) $D^{\alpha }D_{\alpha }\xi _{\beta }  = 0;$
\par iii) $D_{\alpha }D^{\alpha }\xi _{\beta } + k^*_{\beta
\bar{\sigma }}\xi
^{\bar{\sigma }}  = 0.$
\end{pro}
On every compact balanced Hermitian manifold the Ricci 
formula
(\ref{7*}) leads to the
following integral formula
\begin{equation}\label{4.6}
\int_{{\bf M}}  \Vert D_{\alpha }\xi _{\beta } \Vert^2 \,dV
= \int_{{\bf M}}  \Vert
D_{\alpha }\xi _{\bar{\beta }} \Vert^2 \,dV - \int_{{\bf M}}
k^*(\omega
^{\#},\omega ^{\#}) \,dV.
\end{equation}
DEFINITION. A real vector field $\xi $ on a Hermitian
manifold is said to be {\it
affine Hermitian } if it is affine vector field with respect
to the Chern connection
$D$, i.e $L_{\xi }D = 0$. If the linear connection $L_{\xi
}D$ preserves the
complex structure $J$, i.e.
\begin{equation}\label{8}
(L_{\xi }D) \circ J = J \circ (L_{\xi }D) ,
\end{equation}
then we call $\xi $ a {\it complex Hermitian} vector field. \par
In local holomorphic coordinates the condition (\ref{8}) is
equivalent to the
equations
$$
(L_{\xi }D)^{\lambda}_{\alpha \bar{\beta }} = 0; \qquad
(L_{\xi
}D)^{\bar{\lambda }}_{\alpha \beta } = 0.
$$
Using the general formulas expressing $L_{\xi }D$ with the
torsion and
curvature (see \cite{Ya1}) we find
$$
(L_{\xi }D)^{\lambda }_{\alpha \bar{\beta }} = D_{\alpha
}D_{\bar{\beta }}\xi
^{\lambda }; \qquad (L_{\xi }D)^{\bar{\lambda }}_{\alpha
\beta } = D_{\alpha
}D_{\beta }\xi ^{\bar{\lambda }}.
$$
Applying Proposition \ref{pro32} we obtain
\begin{th}\label{th33}
Let $\xi $ be a real vector field on a compact balanced
Hermitian manifold. The
following conditions are equivalent:
\par i) $\xi $ is analytic;
\par ii) $\xi $ is complex Hermitian.
\end{th}
In the Kaehler case  the Chern connection coincides with the
Levi-Civita
connection . From Theorem \ref{th33} we have
\begin{co}\label{co34}
A real vector field $\xi $ on a compact Kaehler manifold is
analytic iff $L_{\xi
}\nabla $ preserves the complex structure.
\end{co}
Every affine Hermitian vector field is complex Hermitian.
>From Theorem
\ref{th33} we obtain
\begin{th}
Every affine Hermitian vector field on a compact balanced
Hermitian manifold is
analytic.
\end{th}
This result extends the well known result that every affine
vector field on a
compact Kaehler manifold is Killing and hence, it is
analytic. \par
We recall that a real vector field $\xi $ is said to be a
Killing vector field if
$L_{\xi }g = 0$. This condition is equivalent to
\begin{equation}\label{9}
(\nabla _X\omega _{\xi })Y + (\nabla _Y\omega _{\xi })X = 0,
X,Y \in {\bf XM}.
\end{equation}
The condition (\ref{9}) implies $\delta \omega _{\xi } = 0$.
\par
Let $\xi $ be a real vector field. We consider the following
real 1-form
$$
\phi = \xi ^{\beta }D_{\beta }\xi _{\alpha }dz^{\alpha } +
\xi ^{\bar{\beta
}}D_{\bar{\beta }}\xi _{\bar{\alpha }}dz^{\bar{\alpha }}.
$$
Using (\ref{4}) we find
\begin{equation}\label{10}
\delta \phi = - 2Re \left[ D_{\beta }\xi _{\alpha }
D^{\alpha }\xi ^{\beta
}\right] -
2Re \left[ \xi ^{\beta }D^{\alpha }D_{\beta }\xi _{\alpha }
+ \xi ^{\beta
}\theta ^{\alpha }D_{\beta }\xi _{\alpha }\right].
\end{equation}
Using the Ricci identities and (\ref{4}) we calculate
\begin{equation}\label{11}
2Re\left[\xi ^{\beta }D^{\alpha }D_{\beta }\xi _{\alpha
}\right] = s(\xi
,\xi ) -
\frac{1}{2}\xi \delta \omega _{\xi } - \frac{1}{2}J\xi
\delta \omega _{J\xi } -
\frac{1}{2}\xi \theta (\xi ).
\end{equation}
On a compact Hermitian manifold we derive from (\ref{11})
and (\ref{10}) the
folowing  formula
$$
\int_{{\bf M}} \left\{ 2Re\left[D_{\beta }\xi _{\alpha
}D^{\alpha }\xi {\beta
} \right] +
s(\xi ,\xi ) - \frac{1}{2}(\delta \omega _{\xi })^2 -
\frac{1}{2}(\delta
\omega
_{J\xi })^2\right\} \,dV -$$ $$- \int_{{\bf M}}\left\{
\frac{1}{2}\xi \theta
(\xi ) + 2Re\left[ \xi
^{\beta }\theta ^{\alpha }D_{\beta }\xi _{\alpha
}\right]\right\} \,dV = 0.
$$
Thus, we proved
\begin{pro}\label{pro36}
If $\xi $ is a real vector field on a compact balanced
Hermitian manifold, then
\begin{equation}\label{12}
\int_{{\bf M}} \left\{ 2Re\left[D_{\beta }\xi _{\alpha
}D^{\alpha }\xi {\beta
} \right] +
k(\xi ,\xi ) - \frac{1}{2}(\delta \omega _{\xi })^2 -
\frac{1}{2}(\delta
\omega _{J\xi })^2 \right\} \,dV = 0.
\end{equation}
\end{pro}
Now, let $\xi $ be a Killing vector field. Using (\ref{3})
from (\ref{9}) we get
\begin{equation}\label{13}
D_{\alpha }\xi _{\beta } + D_{\beta }\xi _{\alpha } = 0.
\end{equation}
Using (\ref{13}) from Proposition \ref{pro36} we get
\begin{pro}\label{pro37}
If $\xi $ is a Killing vector field on a compact balanced
Hermitian manifold, then
\begin{equation}\label{14}
\int_{{\bf M}} \left\{ \Vert D_{\beta }\xi _{\alpha }
\Vert^2 - k(\xi ,\xi )
+ \frac{1}{2}(\delta \omega _{J\xi })^2 \right\} \,dV = 0.
\end{equation}
\end{pro}
As a corollary we obtain the following theorem of Bochner
type
\begin{th}\label{th39}
Let $({\bf M},g,J)$ be a compact balanced Hermitian
manifold.
\par i) If the Chern form $\kappa $ is non-positive definite, 
then
every
Killing vector field $\xi $
on {\bf M} is analytic and satisfies the equality $$k(\xi
,\xi) = \delta (\omega _{J\xi }) = 0.$$
\par ii) If the Chern form $\kappa $ is negative definite,
then
there are no Killing vector
fields other than zero, i.e. the group of isometries of ({\bf
M},g,J) is discrete.
\end{th}
\begin{rem}\label{rem1}
If the Chern form $\kappa $ is negative definite then the
first
Chern class is negative.
By the theorem of S.Kobayashi \cite{Ko} it follows that
there are no holomorphic
vector fields and hence there are no Killing vector fields
with respect to any
Kaehler metric on {\bf M}.
\end{rem}
EXAMPLE. Let ${\bf M} = G/\Gamma$ be a compact quotient of 
a
complex Lie
group $G$ with respect to its discrete subgroup $\Gamma $.
It is well known that
every compact Hermitian manifold with flat Chern connection
is isomorphic with
${\bf M}$ endowed with its flat hermitian structure \cite{Go}. 
If
the group $G$ is non-abelian, then ${\bf M}$ is a non-Kaehler 
balanced
Hermitian manifold. Since $k = k^* = 0$, then every Killing 
vector field is
analytic and by (\ref{4.6}) it is parallel. (There
is a general result of P.Gauduchon \cite{Ga} that on a
compact Hermitian
manifold with negative semi-definite mean curvature $k^*$
every analitic vector
field is parallel (see also \cite{Ko1})).
\begin{rem}\label{rem2}
It is clear from above that Proposition \ref{pro37} and
Theorem \ref{th39} are
valid if we only assume the condition (\ref{13}) which is
weaker than the Killing
condition (\ref{9})
\end{rem}
\section{Harmonic and analytic 1-forms}

A real 1-form $\omega $ is {\it analytic } if
its $(1,0)$-part is holomorphic. In terms of the Chern 
connection this
condition is equivalent
to the condition
\begin{equation}\label{4.1}
D_{\alpha }\omega _{\bar{\beta }} = 0
\end{equation}
We recall that a real 1-form $\omega $ on a Riemannian
manifold is {\it
harmonic} if it is closed and co-closed, i.e.
\begin{equation}\label{4.2}
d\omega = 0; \qquad \delta \omega = 0.
\end{equation}
Using (\ref{3}) the condition $d\omega = 0$ is equivalent to
the following two
conditions:
\begin{equation}\label{4.3}
D_{\alpha }\omega _{\bar{\beta }} - D_{\bar{\beta }}\omega
_{\alpha } = 0.
\end{equation}
\begin{equation}\label{4.4}
D_{\alpha }\omega _{\beta } - D_{\beta }\omega _{\alpha } =
- T_{\alpha \beta
}^{\sigma }\omega _{\sigma }.
\end{equation}
The condition (\ref{4.3}) implies
\begin{equation}\label{4.5}
\delta (J\omega ) = 0.
\end{equation}
We are going to obtain necessary and sufficient conditions
for  a
harmonic
1- form to be analytic and for a holomorphic 1-form to be
harmonic on a compact  balanced Hermitian manifold.For this
purpose we need some
integral
formulas.
\begin{pro}\label{pro41}
For every real 1-form $\omega $ on a compact balanced
Hermitian manifold the
following integral formulas are valid:
\begin{equation}\label{4.7}
\int_{{\bf M}} \frac{1}{2} \Vert D_{\alpha }\omega _{\beta }
- D_{\beta }\omega
_{\alpha } + T_{\alpha \beta }^{\sigma }\omega _{\sigma }
\Vert^2 \,dV  =
\end{equation} $$\int_{{\bf
M}} \left[\Vert D_{\alpha }\omega _{\bar{\beta }} \Vert^2 +
k(\omega
^{\#},\omega
^{\#}) - k^*(\omega ^{\#},\omega ^{\#}) + \frac{1}{2}
t(\omega ^{\#},\omega
^{\#}) \right] \,dV - $$
$$
- \int_{{\bf M}} \left\{ \frac{1}{2} (\delta \omega )^2 +
\frac{1}{2} (\delta (J\omega
))^2 - 2Re \left[ T_{\alpha \beta }^{\sigma }\omega _{\sigma
}(D^{\alpha }\omega
^{\beta } - D^{\beta }\omega ^{\alpha })\right]\right\} \,dV
= 0.
$$
\begin{equation}\label{4.8}
\int_{{\bf M}} \left\{ k(\omega ^{\#},\omega ^{\#}) -
k^*(\omega ^{\#},\omega
^{\#})
+ \frac{1}{2} 2Re \left[ T_{\alpha \beta }^{\sigma }\omega
_{\sigma }(D^{\alpha
}\omega ^{\beta } - D^{\beta }\omega ^{\alpha })\right]
\right\} \,dV +
\end{equation}
$$
 + \int_{{\bf M}} 2Re \left[ T_{\alpha \beta }^{\sigma
}D^{\alpha }\omega _{\sigma
}\omega ^{\beta } \right] \,dV = 0 $$
\end{pro}
\indent {\it Proof:} The formula (\ref{4.7}) follows
immediately from (\ref{12})
and(\ref{4.6}).
\par
To prove (\ref{4.8}) we consider the following real 1-form
$$
\psi = T_{\alpha \beta }^{\sigma }\omega _{\sigma }\omega
^{\beta }dz^{\alpha } +
T_{\bar{\alpha }\bar{\beta }}^{\bar{\sigma }}\omega
_{\bar{\sigma }}\omega
^{\bar{\beta }}dz^{\bar{\alpha }}.
$$
Applying (\ref{4})we have
\begin{equation}\label{4.9}
- \delta \psi = 2Re \left[ D^{\alpha }T_{\alpha \beta
}^{\sigma }\omega _{\sigma
}\omega ^{\beta }\right] + \end{equation} $$ + 2Re \left[
T_{\alpha \beta
}^{\sigma }D^{\alpha }\omega
_{\sigma }\omega ^{\beta }\right] + Re \left[ T_{\alpha
\beta }^{\sigma }\omega
_{\sigma }(D^{\alpha }\omega ^{\beta } - D^{\beta }\omega
^{\alpha })\right]
$$
>From the second Bianchi identity we get
\begin{equation}\label{4.10}
2Re \left[ D^{\alpha }T_{\alpha \beta }^{\sigma }\omega
_{\sigma }\omega ^{\beta
}\right] = s(\omega ^{\#},\omega ^{\#})
- k^*(\omega ^{\#},\omega ^{\#}).
\end{equation}
Substituting (\ref{4.10}) into (\ref{4.9}) and integrating
the obtained
equality over {\bf M} we obtain (\ref{4.8}) \hfill {\bf
Q.E.D.}

We define the tensor $H$ by the equality
\begin{equation}\label{4.11}
H(X,Y) := k(X,Y) - k^*(X,Y) - \frac{1}{2}t(X,Y), \qquad X,Y
\in {\bf XM}
\end{equation}
>From this definition it follows that the tensor $H$ is
symmetric and $J$-invariant.\\
We have
\begin{th}\label{th42}
Let $({\bf M},J,g)$ be a compact balanced Hermitian
manifold.\\[1mm]
\indent i) A harmonic 1-form $\omega $ is analytic iff
$\int_{{\bf M}}H(\omega
^{\#},\omega ^{\#}) \,dV = 0$;\\[2mm]
\indent ii) An analytic 1-form $\omega $ is harmonic iff
$\int_{{\bf M}}H(\omega ^{\#},\omega ^{\#}) \,dV = 0$;
\end{th}

\indent {\it Proof:} The theorem follows from the following
two lemmas:
\begin{lem}\label{lm43}
Let $\omega $ be a harmonic 1-form on a compact balanced
Hermitian
manifold. Then we have
\begin{equation}\label{4.12}
\int_{{\bf M}} \left[ \Vert D_{\alpha }\omega _{\bar{\beta
}} \Vert^2 + H(\omega
^{\#},\omega ^{\#}) \right] \,dV = 0.
\end{equation}
\end{lem}
\begin{lem}\label{lm44}
Let $\omega $ be an analytic 1-form on a compact balanced
Hermitian manifold.
Then we have
$$
\int_{{\bf M}} \left[ \frac{1}{2} \Vert D_{\alpha }\omega
_{\beta } - D_{\beta
}\omega _{\alpha } + T_{\alpha \beta }^{\sigma }\omega
_{\sigma } \Vert^2 +
H(\omega ^{\#},\omega ^{\#})\right] \,dV = 0.
$$
\end{lem}
The proof of Lemma \ref{lm43} follows after substitution of
(\ref{4.3}), (\ref{4.4})
and (\ref{4.5}) into (\ref{4.7}). Combining (\ref{4.8}) with
(\ref{4.7}) and using
(\ref{4.1}) we get the proof of Lemma \ref{lm44}. This
completes the proof of
Theorem \ref{th42}. \hfill {\bf Q.E.D.}
\begin{rem}\label{rem41}
On a Kaehler manifold the tensor $H$ vanishes identically
and Theorem
\ref{th42} implies the well known result that on a compact
Kaehler manifold every
harmonic 1-form is analytic and vice versa.
\end{rem}
Using Hodge theory (see e.g. \cite{Be}) from Theorem
\ref{th42} we obtain
\begin{co}\label{co45}
For a compact balanced Hermitian manifold with zero tensor
$H$ the de Rham
cohomology group $H^1({\bf M},{\bf R})$ is isomorphic to the
Dolbeault
cohomology \\ group $H^{1,0}({\bf M},{\bf C})$.
\end{co}
Now we can state our main result (vanishing theorem of
Bochner type)
\begin{th}\label{th46}
Let $({\bf M},J,g)$ be a compact balanced Hermitian
manifold.
\par i) If the tensor $H$ is positive semi-definite then
every analitic $1$-form
$\omega $ is harmonic and vice versa, every harmonic $1$-
form $\omega $ is
analytic; moreover in the two cases $H(\omega ^{\#},\omega
^{\#}) = 0$.
\par ii) If the tensor $H$ is positive definite on {\bf M}
then:
\par - there are no harmonic $1$-forms on {\bf M};
\par - there are no analytic $1$-forms on {\bf M}.
\end{th}

\indent {\it Proof:} The proof of this theorem follows
immediately from Lemma
\ref{lm43} and Lemma \ref{lm44}. \hfill {\bf Q.E.D.}\\
Applying the Hodge theory we get
\begin{th}\label{th47}
Let $({\bf M},J,g)$ be a compact balanced Hermitian
manifold with positive
definite tensor $H$. Then
\par i) the first Betti number $b_1 = dimH^1({\bf M},{\bf
R}) = 0$;
\par ii) the Hodge number $h^{1,0} = dimH^{1,0}({\bf M},{\bf
C}) = 0$.
\end{th}
Using (\ref{4.6}) we can state Lemma \ref{lm43} as
\begin{lem}\label{l43'}
Let $\omega $ be a harmonic $1$-form on a compact balanced
Hermitian
manifold. Then the following formula is true
\begin{equation}\label{4.13}
\int_{{\bf M}} \left[ \Vert D_{\alpha }\omega _{\beta }
\Vert^2 + k(\omega
^{\#},\omega ^{\#}) - \frac{1}{2}t(\omega ^{\#},\omega
^{\#})\right] \,dV = 0
\end{equation}
\end{lem}
>From (\ref{4.13}) we get the following vanishing theorem of
Bochner type
\begin{th}\label{th48}
Let $({\bf M},J,g)$ be a compact balanced Hermitian
manifold.
\par i) If the tensor $k - \frac{1}{2}t$ is positive semi-definite  
then the
vector field $\omega ^{\#}$ corresponding to any harmonic 
$1$-form
$\omega $ is
holomorphic
vector field and $$(k - \frac{1}{2}t)(\omega ^{\#},\omega
^{\#}) = 0.$$
\par ii) If the tensor $k - \frac{1}{2}t$ is positive
definite on {\bf M} then there are
no harmonic $1$-forms on {\bf M} and consequently $b_1 = 0$.
\end{th}
\begin{rem}
If the tensor $k - \frac{1}{2}t$ is positive definite on
{\bf M} then the first Chern
form $\kappa $ is positive, by the properties of $t$. From
the
Calabi-Yau-Aubin theory
(see\cite{Be}) there exists a Kaehler metric on {\bf M} with
positive Ricci tensor
and the last conclusion in ii) of Theorem \ref{th48} follows
from the classical
vanishing theorem of Bochner.
\end{rem}
It is well known that on a compact Hermitian manifold if the
mean curvature
$k^*$ is positive semi-definite and it is positive definite
in one point  then there
are no holomorphic $(p,0)$-forms (\cite{KW,Ga}, see also
\cite{Wu,Ko1}), in
particular
there are no holomorphic $(1,0)$-forms (the latter fact for
balanced manifolds
follows also from (\ref{4.6})). In view of Theorem
\ref{th46},  we can generalize
the latter fact for compact balanced manifolds as follows
\begin{th}\label{th49}
If on a compact balanced Hermitian manifold the mean
curvature satisfies the
condition
$$k^*(X,X) < k(X,X) - \frac{1}{2}t(X,X), \qquad X \in {\bf
XM}$$
then there are neither holomorphic $(1,0)$-forms nor
harmonic $1$-forms on
{\bf M}.
\end{th}
As a corollary we also have
\begin{th}\label{th10}
Let $({\bf M},J,g)$ be a compact balanced Hermitian
manifold. If the tensors
$H$ and $k^*$ are positive semi-definite on {\bf M} and
$k^*$ is positive definite
in one point then there are no harmonic $1$-forms on {\bf M}
and consequently
$b_1 = 0$
\end{th}

\noindent {\bf Author's address:}\\[2mm]
Georgi Ganchev, \\
Bulgarian Academy of Science, Institute of Mathematics,
Acad. G.Bonchev
Str.,blok 8, 1113 Sofia BULGARIA\\[2mm]
Stefan Ivanov\\
University of Sofia, Faculty of Mathematics and Informatics,
Department of
Geometry,  bul. James Bouchier 5, 1126 Sofia, BULGARIA.\\
E-mail: ivanovsp@fmi.uni-sofia.bg

\begin{thebibliography}{99}
\bibitem{A_B1} L.Alessandrini and G.Bassaneli, {\it Positive
$\partial\bar{\partial}$-closed currents and non-Kaehler
geometry},  to appear
in The J. of Geom. Analysis.
\bibitem{A_B2} L.Alessandrini and G.Bassaneli, {\it Smooth
proper modification
of compact Kaehler manifolds}, Proc. of the Intern. Workshop
on Complex
Analysis, Wupertal $1990$, Complex Analysis, Aspects of
Math. E vol. ${\bf
17}$, Germany $1991$.
\bibitem{A_B3} L.Alessandrini and G.Bassaneli, {\it Metric
properties of
manifolds bimeromorphic to compact Kaehler manifolds}, J.
Diff. Geometry ${\bf
37} (1993), 95-121$.
\bibitem{AGI} B.Alexandrov, G.Grantcharov, S.Ivanov, {\it
Curvature properties
of twistor spaces of quaternionic Kaehler manifolds}, to
appear.
\bibitem{A-Z} H.Akbahr-Zadeh, {\it Transformations
holomorphiquement projectives
des varietes hermitiennes et kaehleriennes}, J. Math. pures
at appl. ${\bf 67}
(1988), 237-261$.
\bibitem{Ba} A.Balas , {\it Compact Hermitian manifolds
of constant holomorphic sectional curvature.}, Math. Z.
${\bf 189}
(1985),193-210.$
\bibitem{Be} A.Besse, {\it Einstein manifolds}, Springer,
Berlin $1987$.
\bibitem{Ga} P.Gauduchon, {\it Fibres hermitiennes a
endomorphisme de
Ricci non-negatif.}, Bull. Soc. Math. France {\bf 105}
$(1977), 113-140.$
\bibitem{Ga1} P.Gauduchon, {\it Structures de Weyl et
theoremes d'anulation sur une variete conforme autoduale},
Ann. Scuola Norm. Sup. Pisa, Serie IV, vol. XVIII, Fasc.
${\bf 4} (1991), 563-629.$
\bibitem{Go} S.Goldberg, {\it Curvature and homology}, New-
York: Academic Press $1962$.
\bibitem{H} F.Hirzebruch, {\it New geometric methods in
algebraic geometry},
Springer-Verlag, Berlin, Heidelberg, New-York $(1966)$.
\bibitem{Ko} S.Kobayashi, {\it Transformation groups in
Differential Geometry},
Springer-Verlag, Berlin, Heidelberg, New-York $(1972)$.
\bibitem{Ko1} S.Kobayashi, {\it Differential Geometry of
complex vector
bundles}, Iwanami Shoten, Publishers and Princeton Univ.
Press. $1987$.
\bibitem{KW} S.Kobayashy and H.Wu, {\it On holomorphic
sections of certain
Hermitian vector bundles}, Math. Ann. {\bf 189} $(1970), 1-
4.$
\bibitem{Mi} M.L.Michelson, {\it On the existence of special
metrics in complex
geometry}, Acta Math. ${\bf 143} (1983), 261-295$.
\bibitem{P} M.Pontecorvo, {\it Complex structures on
quaternionic manifold}, Diff. Geom. and its Appl. ${\bf 4}
(1994), 163-177.$
\bibitem{Wu} H.Wu, {\it Bochner technique in Differential
Geometry}, Math.
Reports, vol.$3$, part $2$, $1988$.
\bibitem{Ya1} K.Yano, {\it The theory of Lie derivatives and
its applications},
Amsterdam, North-Holland $1957$.
\bibitem{Ya2} K.Yano, {\it Differential geometry on Complex
and Almost Complex
spaces}, Pergamon Press, $1965$.
\end{thebibliography}
\end{document}